\documentclass[12pt,a4paper]{amsart}
\usepackage{amssymb,amsfonts,amsthm,amsmath}

\theoremstyle{plain}
\newtheorem{theorem}{Theorem}
\newtheorem{lemma}{Lemma}

\numberwithin{equation}{section} \numberwithin{theorem}{section}
\numberwithin{lemma}{section} \numberwithin{definition}{section}
\numberwithin{corollary}{section}
\numberwithin{proposition}{section} \textheight =24cm
\begin{document}
\title[Gaussian integer solutions]{Gaussian integer solutions for the fifth power taxicab number problem}
\author{Geoffrey B Campbell}
\address{Mathematical Sciences Institute,
         The Australian National University,
         Canberra, ACT, 0200, Australia}

\email{Geoffrey.Campbell@anu.edu.au}

\author{Alexander Zujev}
\address{Department of Physics,
         University of California,
         Davis, CA, USA}

\email{azujev@ucdavis.edu}


\begin{abstract}
The famous open problem of finding positive integer solutions to $a^5 + b^5 = c^5 + d^5$ is considered,
and related solutions are found in two distinct settings: firstly, where $a$ and $b$ are both positive integers with $c$ and $d$ both  Gaussian integers; secondly, where all of $a$, $b$, $c$, and $d$ are Gaussian integers.
\end{abstract}

\maketitle

\section{Introduction} \label{S:intro}

A famous open question \cite{Dickson1999} is the solution in positive integers of
\begin{equation}\label{E:1.1}
w^5 + x^5 = y^5 + z^5.
\end{equation}
While not resolving this for integers generally, we give an infinite set of solutions with integers $w$ and $x$ where the right side $y$ and $z$ are Gaussian integers. Also, we give an infinite set of solutions where all of $w$, $x$, $y$, and $z$ are Gaussian integers. The solutions to  \begin{equation}
w^n + x^n = y^n + z^n,
\end{equation}
\noindent
for $n=4$ are well known and date back to Euler. See Hardy and Wright \cite{Hardy1954} and the ebook by Piezas \cite{Piezas2015}, for example. The case where $n=3$ is solved by the well known and celebrated "Taxicab numbers" named after the famous Hardy and Ramanujan anecdote. (See Hardy \cite{Hardy1927}.)

\section{Solutions where $a$ and $b$ are both positive integers with $c$ and $d$ both Gaussian integers} \label{S:first_thm}

Our first result is encapsulated in the
\begin{theorem} \label{Thm1}
  If the Pell number sequence is 0, 1, 2, 5, 12, 29, 70, 169, 408, 985, 2378,...; then an infinite sequence of solutions to $(\ref{E:1.1})$ is given by:
  \begin{equation}
    3^5 + 1^5 = (2 + i 3)^5 + (2 - i 3)^5,
    \end{equation}
    \begin{equation}
    13^5 + 11^5 = (12 + i 17)^5 + (12 - i 17)^5,
    \end{equation}
    \begin{equation}
    71^5 + 69^5 = (70 + i 99)^5 + (70 - i 99)^5,
    \end{equation}
    \begin{equation}
    409^5 + 407^5 = (408 + i 577)^5 + (408 - i 577)^5,
    \end{equation}
    \begin{equation}
    2379^5 + 2377^5  = (2378 + i 3363)^5 + (2378 - i 3363)^5,
    \end{equation}
    and so on where for $P_n$ the $n$th Pell number, the $n$th equation is
    \begin{equation}
    (P_{2n+3}+1)^5 + (P_{2n+3}-1)^5 = (P_{2n+3} + i (P_{2n+3}+P_{2n+2}))^5 + (P_{2n+3} - i (P_{2n+3}+P_{2n+2}))^5.
    \end{equation}
\end{theorem}
It does seem interesting that the ancient Pell number sequence should figure so neatly in the above set of solutions, with integers on the left, Gaussian integers on the right. The proof of Theorem \ref{Thm1} is by simple expansion.

\section{Solutions where all of $a$, $b$, $c$ and $d$ are Gaussian integers} \label{S:second_thm}

Our second result requires the following identity,
\begin{lemma} \label{L:3.1}
For all real values of $a$, $b$, $c$,
\begin{equation}
(a+b+ ic)^5 + (a-b- ic)^5 - (a-b+ ic)^5 - (a+b- ic)^5 = 80abc (a^2 + b^2 - c^2).
\end{equation}
\end{lemma}
We see  that every Pythagorean triple $a$, $b$, $c$ yields a zero on the right side of the lemma, and hence a Gaussian integer solution of $A^5 + B^5 = C^5 + D^5$. This proves the
\begin{theorem} \label{Thm2}
  Every Pythagorean triple $a$, $b$, $c$ implies a Gaussian integer solution of $A^5 + B^5 = C^5 + D^5$.
  \end{theorem}
Some examples with primitive triples are: \\
$(4,3,5)$ leads to $(7 + 5i)^5 + (1 - 5i)^5$ = $(7 - 5i)^5 + (1 + 5i)^5$, \\
$(12,5,13)$ leads to $(17 + 13i)^5 + (7 - 13i)^5$ = $(7 + 13i)^5 + (17 - 13i)^5$, \\
$(15,8,17)$ leads to $(23 + 17i)^5 + (7 - 17i)^5$ = $(7 + 17i)^5 + (23 - 17i)^5$. \\


\begin{thebibliography}{99}
\bibitem{Apostol1976}
APOSTOL, T.,  "Introduction to Analytic Number Theory",
Springer-Verlag, New York, 1976.
\bibitem{Hardy1927}
HARDY, G. H., "Ramanujan", Cambridge University Press, 1927.
\bibitem{Dickson1999}
DICKSON, L. E., "History of the Theory of Numbers", Vol II, Ch XXII, page 644, originally published 1919 by Carnegie Inst of Washington, reprinted by The American Mathematical Society 1999.
\bibitem{Hardy1954}
HARDY G. H., and WRIGHT E.M., "An Introduction to the Theory of Numbers", 3rd ed., Oxford University Press, London and NY, 1954, Thm. 412.
\bibitem{Piezas2015}
PIEZAS III, T., ebook: https://sites.google.com/site/tpiezas/020 (fifth powers)
\bibitem{Weisstein2015}
WEISSTEIN, Eric W.,"Pell Number", from MathWorld-A Wolfram Web Resource. http://mathworld.wolfram.com/PellNumber.html
\end{thebibliography}
\end{document}